\documentclass[]{amsart}

\usepackage[utf8]{inputenc}
\usepackage{amssymb}
\usepackage{amsmath}
\usepackage{amscd}
\usepackage{epic}
\usepackage{xypic}

        \def\Z{\mathbb{Z}}

        \def\F{\mathbb{F}}
        \def\P{\mathbb{P}}

    \newtheorem{Theo}{Theorem}[section]
	\newtheorem{Defi}[Theo]{Definition}
	\newtheorem{Propo}[Theo]{Proposition}
	\newtheorem{Coro}[Theo]{Corollary}
	\newtheorem{Lemme}[Theo]{Lemma}
	\newtheorem{ex}[Theo]{Example}
    \newtheorem{Rem}[Theo]{Remark}

\newcommand{\biindice}[3]%
 {
 
 {#1}_{\begin{array}[t]{c}
         {\scriptstyle #2} \\
         {\scriptstyle #3}
       \end{array}}
 
 }
\newcommand{\Deg }{\mathrm{Deg\, }}
\newcommand{\im }{\mathrm{Im\,}}
\newcommand{\Ker }{\mathrm{Ker\,}}
\newcommand{\norm }{\mathrm{norm\,}}
\newcommand{\tr }{\mathrm{trace\,}}
\newcommand{\Gal }{\mathrm{Gal\,}}
\newcommand{\ord }{\mathrm{ord\, }} 
\newenvironment{preuve}{\noindent {{\bf Proof}.}}{\hspace*{\fill}$\square$\vskip 8pt}

\author[S. Ballet]{St\'ephane Ballet}
\address{Aix Marseille Univ, CNRS, Centrale Marseille, I2M, Marseille, France \\
Institut de Math\'ematiques de Marseille\\
Case 907, 13288 Marseille cedex 9.}
\email{stephane.ballet@univ-amu.fr}

\author[R. Rolland]{Robert Rolland}
\address{Aix Marseille Univ, CNRS, Centrale Marseille, I2M, Marseille, France \\
Institut de Math\'ematiques de Marseille\\
Case 907, 13288 Marseille cedex 9.}
\email{robert.rolland@acrypta.fr}

\title[On the Weil descent of bi-cyclic Artin-Schreier Function Fields]{On the Weil descent of Artin-Schreier algebraic function fields over finite fields}

         \def\bib#1#2#3#4{\bibitem{#1}{\it #2}{\ #3}{\ #4}}
\date{\today}
\begin{document}

\subjclass{11G20, 14H25}
\keywords{Algebraic function field, finite field, descent of definition field.}

\begin{abstract}
Let us consider a generalized Artin-Schreier algebraic function field extension $F$ of the rational function field $\F_{p^n}(x)$
defined over the finite field extension $K=\F_{p^n}$ of the prime field $\F_p$. We assume that $K$ is algebraically closed in $F$. 
We give general results on the descent over the fields $k= \F_{p^t}$ for $t$ dividing $n$.
Then, we completely handle the bi-cyclic case of the descent over the fields $k_1=\F_{p}$ and $k_2= \F_{p^2}$ of all the sub-extensions of $F$ defined over $\F_{p^4}$. We give explicit examples with small prime numbers $p$.
\end{abstract}

\maketitle


\section{Introduction}\label{un}

\subsection{General context}

The topic of this paper concerns the arithmetic in algebraic function fields in one variable defined over finite fields, particularly in Artin-Schreier-type algebraic function fields.  
We are interested in a generalized Artin-Schreier $F=F_0(y)$ which is defined by an equation $y^{p^s}+y=f(x)$ 
defined over the finite field extension $K=\F_{p^{2s}}$ of the prime field $\F_p$, where $f$ is an element in the rational field $F_0=\F_{p^{2s}}(x)$, 
which satisfies the following conditions (see \cite[Proposition 3.7.10]{stich}):  
\begin{itemize}
\item For each place $P\in F_0$, there is an element  $z\in F_0$ such that 
\begin{equation} \label{AScase1}
v_P(f(x)-(z^{p^s}+z))\geq 0
\end{equation}
or 
\begin{equation}\label{AScase2}
v_P(f(x)-(z^{p^s}+z))=-m \hbox{ with } m>0 \hbox{ and } m\not \equiv 0 \hbox{ mod } p
\end{equation}
where $v_P$ denotes the discrete valuation associated to the place $P$.

\item There exists at least one place $Q\in \P_{F_0}$ with $m_Q>0 $ where $m_Q:=m$. 
\end{itemize}
Note that in this case, $F/F_0$ is Galois with $[F:F_0]=p^s$ and the Galois group of $F/F_0$ is isomorphic to the additive group $G=Gal(F/F_0)= (\Z/p\Z)^s$, 
and $K$ is the full constant field of $F$ i.e. $K$ is algebraically closed in $F$. 

\medskip

We study the descent of the definition field of the sub-extensions $E$ such that $F_0\subseteq E \subseteq F$.
We are interested by the problem of the descent of the definition field which we recall the general problematic.
Let ${\mathcal C}$ be an absolutely irreducible algebraic curve.
We assume that ${\mathcal C}$
is defined over $K$. That is, the ideal of ${\mathcal C}$
can be generated by a family of polynomials 
which coefficients are in $K$.
Let $F=K({\mathcal C})$ be the 
algebraic function field of one variable associated to
${\mathcal C}$. We remark that 
$K$ is the full constant field of $F$.
The descent problem of the definition field of 
the curve ${\mathcal C}$ from $K$ to a field $k\subset K$ is the following: is it possible
to find an absolutely irreducible curve
${\mathcal C}'$ defined over $k$ and birationally equivalent to ${\mathcal C}$?
In other words, is it possible to find a model ${\mathcal C}'$ of the curve ${\mathcal C}$ with coefficients in $k$.
Let us set $F'=k({\mathcal C}')$. Then
$k$ is the full constant field of $F'$.
Algebraically speaking, the descent problem of the definition field of the algebraic function field $F$ from $K$ to $k\subset K$ is the following: is it possible
to find an algebraic function field $F'/k$ defined over $k$ where $k$ is algebraically closed in $F'$ such that $F=F'K$ is the composite field of $F'$ and $K$. 
More precisely, the descent problem 
can be expressed in the following way:

\begin{Defi} Let $F/K$ be an algebraic function field of 
one variable with full constant field $K$.
We consider $F$ as a finite extension of $K(x)$ 
for some $x\in F$  transcendental over $K$. 
We say that the descent of the definition field of $F$ from $K$ to $k$
is possible,
if there exists an extension $F'$ of $k(x)$,
with full constant field $k$, such that
$F$ is isomorphic as $k$-algebra to $F'\otimes_k K$. We will say also
that we can reduce $F/K$ to $F'/k$.
\end{Defi}

\subsection{Organisation}

In Section \ref{deux} we recall certain elementary conditions for the descent to be possible. We also give sufficient conditions for the descent of the sub-extensions of a Galois extension to be possible. This part comes from the descent theory, developed by A. Weil (cf. \cite{weil}) and by J.-P. Serre (cf. \cite{serre}). In Section 3, we specialize the descent problem to case of generalized Artin-Schreier extensions. In particular, we mainly give an effective necessary and sufficient condition for descent to be possible in this case (cf.Theorem \ref{synth}). In Section 4, we study the generalized Artin-Schreier extensions defined by the additive polynomial $T^{p^s}+T$ with $s=2$ and $p$ an arbitrary prime number.

\section{General results on the descent of the definition field}\label{deux}

From now on, $k$ is a perfect field and $K$ is a finite Galois extension of $k$ such that 
$k\subset K\subset U$, where $U$ is a fixed algebraically closed field. Then the rational function field $K(x)$ over the rational function field $k(x)$ also is Galois. We will denote by $\Gamma=\Gal(K(x)/k(x))$ the Galois group of the rational function field $K(x)$ over the rational function field $k(x)$.
The group $\Gamma$ also is the Galois group of the Galois extension $K/k$ up to isomorphism (cf. Lemma \ref{GaloisCorpsFoncRat}).
Let us consider an algebraic function field $F/K(x)$ which is Galois over $K(x)$.
Let $G=Gal(F/K(x))$ be the Galois group of $F$ over $K(x)$.
For any algebraic function field $F/K$, we assume tacitely 
that $K$ is its full constant field and that $F$ is a subfield of
the algebraically closed field $U$.

\medskip

Now, we are going to study some general necessary and sufficient conditions
allowing the descent of the definition field in our framework. First, let us recall an useful result (cf. \cite[p.10]{sadi}):

\begin{Lemme}\label{GaloisCorpsFoncRat}
Let $K$ be a Galois extension of the perfect field $k$. Then, we have:

\begin{itemize}\label{LemmeDescenteGalois}
\item[(1)] $Gal(K(x)/k(x))\simeq Gal(K/k)$ where $x$ is transcendantal over the field $k$.
\item[(2)] Let $F/K$ be an algebraic function field of one variable with full constant field $K$ and $F'/k$ 
such that $F/K=F'/k\otimes_k K$ i.e. $F/K$ can be reduced over $k$. Then $F/F'$ is Galois and $\Gamma=Gal(F/F')\simeq Gal(K(x)/k(x)) \simeq Gal(K/k)$.
\end{itemize}
\end{Lemme}

Now, let us recall the following sufficient condition (cf. \cite[Theorem 2.2]{ballbrigro}).

\begin{Theo}\label{desc} Let $K$ be a Galois extension of the perfect field $k$ and let $F/K$ be an algebraic function field. 
If the extension $F/k(x)$ is Galois with Galois group ${\mathcal G}$ and 
if ${\mathcal G}$ is the semi-direct product ${\mathcal G}=G \rtimes \Gamma$,
where $G=\Gal(F/K(x))$ and $\Gamma=\Gal(K(x)/k(x))$,
then the descent of the definition field from
$K$ to $k$ is possible. 
\end{Theo}

Now, we focus on the reciprocal problem namely the necessary conditions.

\begin{Theo}\label{rdesc}Let $K$ be a Galois extension of the perfect field $k$ with Galois group $\Gamma$ and let $F$ be a Galois extension of the rational function field $K(x)$ with Galois group $G$. Assume that
it is possible to reduce $F/K$ to $F'/k$. 
Then:

\medskip

(1) $F/k(x)$ is Galois  with Galois group ${\mathcal G}= G\rtimes \Gamma$. 

\medskip

(2) Moreover, if $H$ is a subgroup of $G$ such that the fixed field $E/K$ can be reduce over the algebraic function field $E'/k$ then $F/E'$ is Galois 
with Galois group ${\mathcal H}= H\rtimes \Gamma$.

\medskip

In particular $G$ and $H$ are stabilised under the action by inner automorphism of the group $\Gamma$.
\end{Theo}

\medskip

\begin{preuve}
(1) The extension $F/F'$ is Galois, and its Galois group is isomorphic to $\Gamma$ by Lemma \ref{LemmeDescenteGalois}. 
So that we can define an action $\Psi$ from $\Gamma$ on the Galois group $G=\Gal(F/K(x))$ by inner automorphisms:
$$\Psi(g)=\gamma^{-1}\circ g \circ \gamma,$$
where $\gamma \in \Gamma$ and $g \in G$.

Then, the semi-direct product $G\rtimes \Gamma$ is a group of 
$k(x)$-automorphisms 
of $F$. The cardinality of this group is exactly the degree of 
$F$ over $k(x)$. Hence we can conclude that $G\rtimes \Gamma$ is
the group of automorphisms of $F$ over $k(x)$ and that
the extension $F/k(x)$ is Galois, with Galois group ${\mathcal G}= G\rtimes \Gamma$.
Then we have the short exact sequence
$$1 \rightarrow G \overset{i_1}{\rightarrow} {\mathcal G} 
\overset{\underset{\pi}{\overset{s_1}{\curvearrowleft}}}{\rightarrow} \Gamma \rightarrow 1,$$
where $s_1$ is the corresponding section.

(2) First, it is clear that $F/E$ is Galois with Galois group $H$. As $E/K$ can be reduced over $E'/k$, the extension $E/E'$ is Galois, and its Galois group is isomorphic to $\Gamma$.  So, $F'$ is a finite extension of $E'$ of degree $[E:E']=[F:F']$. So, $F/E'$ is Galois with Galois group ${\mathcal H}=H\rtimes \Gamma$ and we have the short exact sequence: 
$$1 \rightarrow H \overset{i_2}{\rightarrow} {\mathcal H} 
\overset{\underset{\pi}{\overset{s_2}{\curvearrowleft}}}{\rightarrow} \Gamma \rightarrow 1,$$ such that 

$$\begin{array}{ccccccccc}
1 &\rightarrow &G &\overset{i_1}{\rightarrow} &{\mathcal G} &
\overset{\underset{\pi}{\overset{s_1}{\curvearrowleft}}}{\rightarrow} &\Gamma &\rightarrow &1 \\
    &                  &\cup &   & \cup & & & &   \\
1 & \rightarrow &H & \overset{i_2}{\rightarrow} & {\mathcal H} &
\overset{\underset{\pi}{\overset{s_2}{\curvearrowleft}}}{\rightarrow} & \Gamma & \rightarrow &1
\end{array}
$$
where $s_1$ and $s_2$ are the corresponding sections (cf. Figure \ref{fig3}). 

\begin{equation}
\label{fig3}
\xymatrix{
F^\prime \ar@{--}[r]^\Gamma
 &F&F\ar@{-}[dd]^G \\
E' \ar@{--}[u]  \ar@{--}[r] \ar@{--}[ur]^{\mathcal H}&E\ar@{-}[u]_H&\\
k(x) \ar@{--}[u] \ar@{-}[uur]_<<<<<<<<<{\mathcal G} \ar@{-}[r]_\Gamma &K(x)\ar@{-}[u]&K(x)
}
\end{equation}

\end{preuve}

\section{Case of generalized Artin-Schreier extensions}

\subsection{General results}

Let us mention the following well known result (cf. \cite{stich}, Proposition III.7.10) concerning the generalized Artin-Schreier extensions:

\begin{Propo}\label{utilAS} 
Consider an algebraic function field $L/K$ with constant field $K$ of characteristic $p>0$,
and a linearized separable polynomial $A(x)\in K[x]$ of degree $p^s$ which has all its roots in $K$ and $u \in L$. Assume that the polynomial $A(x)-u$  is absolutely irreducible and set $F=L(y)$ with $A(y)=u$.
Then $F/L$ be an elementary abelian $p$-extension of degree $p^s$ and
 the Galois group ${\Gal(F/L)}$ of $F/L$ is defined in the following way : for any $P(y)\in F$, namely a polynomial of degree $<p^s$
 and with coefficients in $L$, and any $\alpha$ in $K$, let us define a function $f_{\alpha}$ from $F$ to $F$
 by $f_{\alpha}(P(y))=P(y+\alpha)$; then
 $${\Gal(F/L)}=\{f_{\alpha} ~|~ \alpha \in K \hbox{ and } A(\alpha)=0\}$$
 where the operation is the composition of functions.
 \end{Propo}

\begin{Coro}
 The Galois group $\Gal(F/L)$ is isomorphic to the additive subgroup $G$ of $K$
 $$G=\{\alpha ~|~ \alpha \in K \hbox{ and } A(\alpha)=0\}.$$
\end{Coro}

\subsection{Case of $A(x)\in k[x]$ and $L=K(x)$}

In this case, $F$ can be reduced to $F'$ over $k$ and the extension $F/F'$ is Galois with Galois group (cf. Lemma \ref{LemmeDescenteGalois})
$$\Gal(F/F')\simeq \Gal(K(x)/k(x)) \simeq \Gal(K/k)=\Gamma.$$
Let us remark that as $A(x)\in k(x)$, for any $\sigma \in \Gamma$ and any $\alpha \in G$
the following holds
$$A\left(\strut\sigma(\alpha)\right)=\sigma\left(\strut A(\alpha)\right)=0.$$
Hence, if $\alpha \in G$ then $\sigma(\alpha)\in G$.

\medskip

Let us consider the following application $\Theta$
from $\Gamma$ in the automorphisms of the group $G$. If $\sigma\in \Gamma$
and $f_{\alpha}\in G$, let us define 
$$\Theta(\sigma)(f_{\alpha})= f_{\sigma(\alpha)},$$
namely
$$\Theta(\sigma)(f_{\alpha})(P)(y)=P(z+\sigma(\alpha)).$$

As $F$ can be reduced, we know that $\Gamma$ acts on the Galois group $G=\Gal(F/K(x))$ by the inner automorphisms
$$\Omega(\sigma)(f_{\alpha})=\sigma^{-1}\circ f_{\alpha}\circ \sigma.$$

\begin{Propo}\label{ident}
For any $\sigma \in \Gamma$ and any $\alpha \in G$, the following holds
$$\Omega(\sigma)=\Theta(\sigma^{-1}).$$
\end{Propo}

\begin{preuve}
 Let $P\in F$, then
 $$P(y)=\sum_{0\leq i <p^s}a_iy^i,$$
$$\sigma(P)(y)=\sum_{0\leq i <p^s} \sigma(a_i)y^i,$$
$$f_{\alpha}\circ\sigma(P)(y)=\sum_{0\leq i <p^s} \sigma(a_i)(y+\alpha)^i,$$
$$\sigma^{-1}\circ f_{\alpha}\circ\sigma(P)(y)=\sum_{0\leq i <p^s} a_i(y+\sigma^{-1}(\alpha))^i,$$
$$\sigma^{-1}\circ f_{\alpha}\circ\sigma(P)(y)=\Theta(\sigma^{-1})(f_{\alpha})(P)(y),$$
$$\sigma^{-1}\circ f_{\alpha}\circ\sigma=\Theta(\sigma^{-1})(f_{\alpha})=f_{\sigma^{-1}(\alpha)}.$$
\end{preuve}

\begin{Theo}\label{synth}
Let $H$ be a subgroup of $G$ and $E$ the subfield of $F$ such that $H=\Gal(F/E)$.
Then $E$ can be reduced over $k$ if and only if for any $\sigma \in \Gamma$, $H$ is globally invariant by
$\Theta(\sigma)$.
\end{Theo}

\begin{preuve}
This is a direct consequence of Proposition \ref{ident} and Theorems \ref{desc}
and \ref{rdesc}.
\end{preuve}

\begin{Coro}\label{stableFrob}
With the notations of Theorem \ref{synth},
if $k$ is the finite field $\mathbb{F}_p$, $K$ the finite field $F_{p^n}$ and $\sigma$
the Frobenius function on $K/k$, the field $E$ can be reduced over $k$ if and only if
the subgroup $H$ of $G$ is globally invariant by
$\Theta(\sigma)$.
\end{Coro}

\section{Generalized Artin-Schreier specific case: $A(T)=T^{p^s}+T$}

Now, let us consider the particular case of generalized Artin-Schreier extensions presented in Introduction. More precisely, we consider the generalized Artin-Schreier extension $F/F_0$ with equation $y^{p^s}+y=f(x)$ defined over the finite field extension $K=\F_{p^{2s}}$ of the prime field $\F_p$, where $f$ is an element in the rational field $F_0=\F_{p^{2s}}(x)$ satisfying the required conditions
(cf. Section \ref{un} Conditions (1) or (2)). In this case, the Galois group $G$ of the extension $F/ K$ over $F_0$ is such that $$G={\Gal(F/F_0)}=\{f_{\alpha} ~|~ \alpha \in K \hbox{ and } \alpha^{{p^s}}+\alpha=0\}.$$

\medskip

By the previous section, we know that for any subgroup $H$ of the Galois group $G$, the fixed field $E_H\subset F/K$ can be reduced over $k_t=\F_{p^t}$ with $t$ dividing $s$ if and only if $H$ can be stabilized over the action of $\Gamma_t=Gal(K/k_t)$, namely under the action of the Frobenius automorphism $\phi_t$ of $\Gamma_t=Gal(\F_{p^{2s}}/k_t)$ defined by: 
 
 $$
 \begin{array}{cccc}
 \phi_t :&  \F_{p^{2s}} & \longrightarrow & \F_{p^{2s}}\\
  & \alpha & \mapsto & \alpha^{p^t} 
 \end{array}
 $$
 where $t$ divides $2s$. Note that  $\phi_t$ fixes the ground field $k_t=\F_{p^{t}}$, on which we are interested to descend.
 
 \subsection{Explicit general description of the group $G$}
 
 In this section, we describe the group $G$ for any integer $s$ and for an odd prime $p$.
 Let $\beta$ be a primitive element of $K=\mathbb{F}_{p^{2s}}$. An element $\beta^u$
is in the group $\mathbf{G}$ if and only if:
\begin{equation}\label{ing}
\beta^{up^s}+\beta^u=0.
 \end{equation}

Then the following theorem holds:

\begin{Theo}\label{Ggroup}
Let $\mathbf{G}$ be the additive group of elements $a$ in $\mathbb{F}_{p^{2s}}$
such that $$a^{p^s}+a=0.$$
The $p^s$ elements of $\mathbf{G}$ are:
$$\left\{\strut 0, \beta^{\frac{1}{2}(p^s+1)}, \beta^{\frac{3}{2}(p^s+1)},\beta^{\frac{5}{2}(p^s+1)},
\cdots, \beta^{\frac{2p^s -3}{2}(p^s+1)}\right\}.$$
\end{Theo}

\begin{preuve}
As $p$ is odd and $\beta$ is a primitive element:
$$-1=\beta^{\frac{1}{2}(p^{2s}-1)}.$$
Then equation (\ref{ing}) can be written:
\begin{equation}\label{ing1}
\beta^{up^s}=\beta^{u+\frac{1}{2}(p^{2s}-1)}.
\end{equation}
The elements $u$ satisfying equation (\ref{ing1}) are the elements $u$ such that:
$$u(p^s-1) \equiv \frac{1}{2}(p^{2s}-1) \mod (p^{2s}-1),$$
$$u \equiv \frac{1}{2}(p^{s}+1) \mod (p^{s}+1).$$
\end{preuve}

Let us determine now the elements of $k_t$ and
the elements not $0$ of $\mathbf{G}$ which are in $k_t$.

\begin{Theo} Let $\beta$ be a primitive element of $K=\mathbb{F}_{p^{2s}}$ and let $t$ be a divisor of $2s$. Then 
$$k_t=\F_{p^t}=\left \{\strut 0, \beta^{\frac{p^s-1}{p^t-1} (p^s+1)}, \beta^{2 \frac{p^s-1}{p^t-1} (p^s+1)}, \cdots,
\beta^{(p^t-1)\frac{p^s-1}{p^t-1} (p^s+1)} \right \}.$$
\end{Theo}

\begin{preuve}
The elements of $k_t$ are
those satisfying:
$$\phi_t(\beta^u)=\beta^u,$$
$$\beta^{p^tu} =\beta^u,$$
$$p^tu\equiv u \mod (p^{2s}-1),$$
$$u(p^t-1) \equiv 0 \mod (p^{2s}-1),$$
$$u \equiv 0 \mod \left (\frac{p^s-1}{p^t-1} (p^s+1) \right).$$
\end{preuve}

\begin{Coro}\label{kt}
The non-zero elements  of $k_t$ are of the form $\beta^{2 n \frac{p^s+1}{2}}$ where $n$ is an integer.
\end{Coro}

\begin{Theo}
$$k_t\cap \mathbf{G}=\{0\}.$$
\end{Theo}

\begin{preuve}
A non-zero element $u\in k_t$ is such that
$$\phi_t(u)=u.$$
and consequently
$$(\phi_t)^{\frac{s}{t}}(u)=u^{p^s}=u.$$
Hence, a non-zero element $u\in k_t$ is not in $G$.
\end{preuve}

\begin{Rem}
The non-zero elements
of $G$ are by Theorem \ref{Ggroup} of the form $\beta^{m \frac{p^s+1}{2}}$ with $m$ odd
and the non-zero elements of $k_t$ are of the form $\beta^{2n \frac{p^s+1}{2}}$.
\end{Rem}

\subsection{The subgroups of $G$}
Let $a_1$ be a non zero element of $G$.
The set $$G_1=\{0,a_1,2a_1,\cdots,(p-1)a_1\}$$ is a subgroup of $G$. Now
let $a_2$ be a non-zero element of $G$ which is not in $G_1$,
the set $$G_2=\{0,a_2,2a_2,\cdots,(p-1)a_2\}$$ is also a subgroup of $G$.
Then, we can define $n=\frac{p^{s}-1}{p-1}$ such subgroups such that:
$$\bigcup_{i=1}^nG_i=G,$$
and if $i\neq j$
$$G_i \cap G_j=\{0\}.$$

\subsection{Study of the case $s=2$, $t=1$}

We study separately the cases where $p=2$ and $p$ is an odd prime.

\subsubsection{Case $p\neq 2$ and $p-1$ not divisible by $4$}

\begin{Theo}
When $p\neq 2$ and $p-1$ is not divisible by $4$, no subgroup $G_i$ is stable by the Frobenius $\phi_1$.
Moreover the set of the $p+1$ subgroups $G_i$ is the union of $\frac{p+1}{2}$ subsets of
two distinct elements $\{G_i,G_j\}$ such that
$$\phi_1(G_i)=G_j \hbox{ and } \phi_1(G_j)=G_i.$$
\end{Theo}

\begin{preuve}
Let $a$ be a non-zero element of $G$ and $G_1$ the additive subgroup
of $\mathbf{G}$ generated by $a$.
First we compute in $\mathbb{F}_{p^4}$.
As $a\in G$ we have
$$a^{p^2}=-a.$$
Let us suppose that $a$ is such that
$a^p=ua$ with $u \in \mathbb{F}_p$, (namely $\phi_1(a) \in G_1$).
Then
$$(ua)^p=a^{p^2}=-a.$$
But
$$(ua)^p=u^p a^p= ua^p=u^2a.$$
Then
$$u^2a=-a,$$
and $u^2=-1$ in $\mathbb{F}_{p^4}$. But by hypothesis $u$ is in the prime field $\mathbb{F}_{p}$
then we also have $u^2=-1$ in $\mathbb{F}_{p}$.
Consequently, il $p-1$ is not divisible by $4$ such a $u$ cannot exist by a consequence of Wilson theorem and
$G_1$ is not stable by the Frobenius $x \rightarrow x^p$.

\medskip

Let $b=\phi_1(a)$. Then $b \in G_2$ where $G_2$ is distinct of $G_1$.
Let $ua$ (where $u\in k_1$) an another element of $G_1$.
$$\phi_1(ua)=(ua)^p=ua^p=ub.$$
Then $\phi_1(ua) \in G_2$ and $(\phi_1)$ is a bijective application from
$G_1$ onto $G_2$. Remark that as $\phi_1^2=-Id$, if $\phi_1(a)=b$ we have $\phi_1(b)=-a$.
\end{preuve}

\subsubsection{Case $p-1$ divisible by $4$}

If $p-1$ is divisible by $4$ there exist two square root of $-1$,
$u_1=\left(\strut(p-1)/2\right)!$ and $-u_1$. If the equation
$x^p=u_1x$ has a solution $a$, it has $p$ solutions which are
all the elements of the additive subgroup generated by $a$:
$$(va)^p=va^p=vu_1a=u_1(va).$$

\begin{Theo}
When $p\neq 2$ and $p-1$ divisible by $4$, among the $p+1$ subgroups $G_i$ of $G$, there is exactly two
subgroups stable by the Frobenius $\phi_1$
(having respectively for elements the solutions of $x^p=u_1x$ and the solutions of $x^p=-u_1x$).
Moreover, if a $G_i$ is not stable by the Frobenius $\phi_1$, its image by $\phi_1$
is another subgroup $G_j$ and the image by $\phi_1$ of $G_j$ is $G_i$.
$\phi_1$ is a bijective application from $G_i$ onto $G_j$.
\end{Theo}

\begin{preuve}
We know that if $\beta$ is a primitive element in $\mathbb{F}_{p^4}$ we have
$$-1=\beta^{\frac{p^4-1}{2}},$$
and as $p^4-1$ is divisible by $4$, $-1$ has two square roots
$$\beta^{\frac{p^4-1}{4}},$$
and
$$\beta^{\frac{p^4-1}{4}+\frac{p^4-1}{2}}$$
in $\mathbb{F}_{p^4}$. These square root are in $\mathbb{F}_p$
if and only if $p-1$ is divisible by $4$.
In this case let us solve the equation
$$a^{p-1}=\beta^{\frac{p^4-1}{4}}.$$
We get
$$a^{p-1}=[\beta^{\frac{(p^2+1)(p+1)}{4}}]^{p-1}.$$
A solution is
$$a=\beta^{\frac{(p^2+1)(p+1)}{4}}. $$
Then the additive subgroup generated by this $a$ is invariant by $\phi$.

An analogous computation gives the second invariant subgroup
which is the additive subgroup generated by
$$b=\beta^{\frac{3(p^2+1)(p+1)}{4}}.$$
\end{preuve}

\subsubsection{Case where $p=2$}
In the case where $p=2$ and $s=2$, we have $\F_{2^{2s}}=\F_{16}$ and the question is: what are the intermediate extensions of $F/\F_{2^{2s}}(x)$ that can be reduced on $\F_4$ and on $\F_2$? First, there exist $\frac{p^s-1}{p-1}=3$ sub-extensions $F_1$, $F_2$ and $F_3$ of degree $p$ of $F/\F_{2^{2s}}$. 
Let us denote $G_1$, $G_2$, and $G_3$ their three respective Galois subgroups of $G=\Gal(F/\F_{2^{2s}})=\F_4\simeq\Z/p\Z\times \Z/p\Z$.  Let us set $\F_4:=\{0,1,\alpha,\alpha^2\}$ where $\alpha$ is the primitive root of the irreducible polynomial $p_{\alpha}(x)=x^2+x+1$ in $\F_2[X]$. Let us set $G_1:=\{0,1\}=\F_2$,  $G_2:=\{0,\alpha\}$ and  $G_3:=\{0,\alpha^2\}$. By the results of Section \ref{stableFrob}, to know if the extensions $F_1/\F_{2^{2s}}$, $F_2/\F_{2^{2s}}$ and $F_3/\F_{2^{2s}}$ can be reduced on $\F_{p^t}$, it is sufficient to see if the respective sub-groups $G_1$, $G_2$, and $G_3$ are stabilized under the action of the Frobenius $\phi_t$ with $t=2$ or $t=1$. It is obvious that  $G_1$, $G_2$, and $G_3$ are stabilized by $\phi_2$ because $G$ is fixed by $\phi_2$. Moreover, we see that $\phi_1(G_2)=G_3$ and $\phi_1(G_3)=G_2$ because 
$\phi_1(\alpha)=\alpha^2$ and $\phi_1(\alpha^2)=(\alpha^2)^2=\alpha$. So, we obtain the following result:

\begin{Theo}
The subgroups $G_1:=\{0,1\}=\F_2$,  $G_2:=\{0,\alpha\}$ and  $G_3:=\{0,\alpha^2\}$ of $Gal(F/\F_{2^{2s}}(x))$ are stabilized under the action of $\phi_2$.
The subgroup $G_1:=\{0,1\}$ is fixed by $\phi_1$; $\phi_1(G_2)=G_3$ and $\phi_1(G_3)=G_2$.
\end{Theo}

\begin{Coro}
All the sub-extensions $F_1/\F_{2^{2s}}(x)$, $F_2/\F_{2^{2s}}(x)$ and $F_3/\F_{2^{2s}}(x)$ of $F/\F_{2^{2s}}(x)$ can be reduce on the definition field $\F_4$. The extension $F_1/\F_{2^{2s}}(x)$ is the unique extension which can be reduce on the definition field $\F_2$.
\end{Coro}

\begin{Rem} The generalized Artin-Schreier extension $F/\F_{2^{2s}}(x)$ defined with the additive polynomial  $A(T)=T^{p^s}+T^p$ can always be also defined on the definition field $\F_{p^t}$ with $t$ dividing $2s$. But in the case where $p=2$, the algebraic function field $F/\F_{p^{t}}(x)$ is Galois if $t=s$ by Theorem \cite[Proposition 3.7.10]{stich}. Hence, we also can deduce that the descent over the definition field $\F_{p^s}$ of all the intermediate extensions of $F/\F_{2^{2s}}(x)$ is possible. In the general case where $t\neq s$ for any prime $p$ or $t=s$ with $p\neq2$, a natural question is to know if $F/\F_{p^{t}}(x)$ is Galois over $\F_{p^{t}}(x)$ when all its sub-extensions of $F/\F_{p^{2s}}(x)$ 
can be reduced over $\F_{p^{t}}$ (all the subgroups of the Galois group $G=Gal(F/\F_{2^{2s}}(x))$ are stabilized by the Frobenius $\phi_t$) with a fixed $t$ dividing $2s$.
\end{Rem}

\appendix

\section{Examples of invariant subgroups when $p-1$ is divisible by $4$}
In the following examples, $P(x)$ is a primitive polynomial of degree $4$ such that
$$\mathbb{F}_{p^4}=\mathbb{F}_{p}(x)/\left (\strut P(x) \right)$$
and $\beta$ a primitive element, root of $P(x)$. Let us denote $u_1$ and $u_2$ the square roots of $-1$ in $\F_p$.

\subsection{Case p=5}
$$P(x)=x^4+x^3+x+3$$
In this case $u_1=2$ and $u_2=3$.
Let us  solve
$x^4=u_1$ and $x^4=u_2$ in $\mathbb{F}_{625}$.
We find the subgroups
$$\{0,  \beta^{585}, \beta^{429}, \beta^{117}, \beta^{273}\}$$
and
$$\{0,  \beta^{39},  \beta^{507},\beta^{195},  \beta^{351}\}.$$

\subsection{Case p=13}
$$P(x)=x^4+x^3+x+2$$
In this case $u_1=5$ and $u_2=8$.
Let us solve
$x^{12}=u_1$ and $x^{12}=u_2$ in $\mathbb{F}_{13^4}$.

We find the subgroups
$$\{0,\beta^{4165},\beta^{6545},\beta^{13685},\beta^{8925},\beta^{25585},\beta^{16065},$$
$$\beta^{1785},\beta^{11305},\beta^{23205},\beta^{27965},\beta^{20825},\beta^{18445}\}$$
and
$$\{0,\beta^{5355},\beta^{7735},\beta^{14875},\beta^{10115},\beta^{26775},\beta^{17255},$$
$$\beta^{2975},\beta^{12495},\beta^{24395},\beta^{595},\beta^{22015},\beta^{19635}\}.$$

\subsection{Case p=17}
$$P(x)=x^4+x^3+5$$
In this case $u_1=4$ and $u_2=13$.

We find the subgroups
$$\{0,\beta^{22185},\beta^{53505},\beta^{6525},\beta^{1305},\beta^{27405},\beta^{37845},\beta^{16965},\beta^{32625},$$
$$\beta^{74385},\beta^{58725},\beta^{79605},\beta^{69165},\beta^{43065},\beta^{48285},\beta^{11745},\beta^{63945}\}$$
and
$$\{0,\beta^{19575},\beta^{50895},\beta^{3915},\beta^{82215},\beta^{24795},\beta^{35235},\beta^{14355},\beta^{30015},$$
$$\beta^{71775},\beta^{56115},\beta^{76995},\beta^{66555},\beta^{40455},\beta^{45675},\beta^{9135},\beta^{61335}\}.$$

The computations were done thanks to the private software SIMBA developped by Robert Rolland and Ren\'e Smadja.

\end{document}